\documentclass[11pt]{article}

\usepackage{fourier}
\usepackage{amsmath}
\usepackage{amssymb}
\usepackage{amsthm}
\usepackage{colonequals}
\usepackage[pdftex, bookmarks, final, colorlinks, linkcolor=black, citecolor=black, urlcolor=blue]{hyperref}
\usepackage{pdfsync}
\usepackage{tensor}
\usepackage[english]{babel}
\usepackage{mathtools}
\usepackage{geometry}
\DeclareMathOperator{\cf}{cf}
\DeclareMathOperator{\dom}{dom}
\DeclareMathOperator{\id}{id}
\DeclareMathOperator{\ran}{ran}

\newcommand{\bigset}[1]{\Bigg\{ #1 \Bigg\}}
\newcommand{\closedSpokes}[2]{\closure{\mathbf{Sp}}(#1, #2)}
\newcommand{\closedSpokesAtInfinity}[1]{\closure{\mathbf{Sp}}^\infty(#1)}
\newcommand{\closure}[1]{\overline{#1}}
\newcommand{\cofinality}[1]{\cf(#1)}
\newcommand{\compactSubsets}[1]{\mathcal{K}(#1)}
\newcommand{\df}{\coloneqq}
\newcommand{\domain}[1]{\dom(#1)}
\newcommand{\family}[3]{(#1)_{#2 \in #3}}
\newcommand{\functions}[2]{{}^{#1} #2}

\newcommand{\locallyContained}[1]{\subseteq_{#1}}
\newcommand{\locallyContainedSup}[2]{\subseteq_{#1}^{#2}}
\newcommand{\nhoodCore}[1]{N_{#1}}
\newcommand{\nhoodCoreSup}[2]{N_{#1}^{#2}}
\newcommand{\nhoodFilter}[1]{\mathcal{N}_{#1}}
\newcommand{\nhoodFilterSup}[2]{\mathcal{N}_{#1}^{#2}}
\newcommand{\range}[1]{\ran(#1)}
\newcommand{\set}[1]{\{ #1 \}}
\newcommand{\spokes}[2]{\mathbf{Sp}(#1, #2)}
\newcommand{\spokesAtInfinity}[1]{\mathbf{Sp}^\infty(#1)}
\newcommand{\Tone}{\text{T}_1}
\newcommand{\transSequence}[3]{(#1)_{#2 < #3}}

\newtheorem{corollary}{Corollary}[section]
\newtheorem{lemma}[corollary]{Lemma}

\newtheorem{question}[corollary]{Question}
\newtheorem{theorem}[corollary]{Theorem}
\theoremstyle{definition}
\newtheorem{definition}[corollary]{Definition}

\title{Convergence properties and compactifications}
\author{Robert Leek \\
University of Oxford}

\begin{document}
\maketitle

\abstract{In this paper, we will use investigate the existence of compactifications with particular convergence properties - pseudoradial, radial, sequential and Fr\'echet-Urysohn - through the use of spoke systems.}
\\\\
\textbf{Keywords}: compactification, Fr\'echet-Urysohn, one-point compactification, pseudoradial, radial, \linebreak
sequential, small cardinals, spoke, spoke system
\\
\textbf{MSC (2010)}: 03E17, 54D35, 54D45, 54D55

\section{Introduction}

In \cite{Leek_internal_2014}, we exhibited a local characterisation of radiality using spoke systems, which are collections of subspaces with neighbourhood bases well-ordered by reverse inclusion at a specified point, that together reconstruct the original neighbourhood filter at that point.
Using this characterisation, we can investigate how to compactify a locally compact space whilst preserving it's convergence properties - the main properties under consideration will be radiality and Fr\'echet-Urysohn, although we will exhibit some results for pseudoradiality and sequentiality too.
We will also demonstrate an alternative characterisation of radiality, using cofinal collections of spokes under local containment.

We shall first investigate one-point compactifications, since if we can compactify and preserve radiality or the Fr\'echet-Urysohn property then the one-point compactification will suffice.
From this, we will see how to improve this result to finite and countable compactifications.
Finally, we use small cardinals to find conditions for uncountable sequential and Fr\'echet-Urysohn compactifications.

For the first two sections of this article, we will not be assuming any separation conditions on our topological spaces.
We recall the following definitions and theorems from \cite{Leek_internal_2014}:

\begin{definition} \hfill
\begin{itemize}
\item A \emph{transfinite sequence} is a net with well-ordered domain, typically indexed by an ordinal with the $\in$-ordering.
If $x$ is a point in a topological space, then a transfinite sequence in $X$ is said to \emph{converge strictly} to a point $x$ in a space if it converges to $x$ and $x$ is not in the closure of any of the proper initial segments of the transfinite sequence.

\item We say that a topological space $X$ is \emph{radial} at a point $x$ if for every subset $A$ of $X$ that contains $x$ in its closure, there is a transfinite sequence converging to $x$ whose range lies in $A$.
If a space is radial everywhere then we call it a \emph{radial} space.

By replacing the transfinite sequences with ordinary $\omega$-indexed sequences above, we arrive at the definition of the \emph{Fr\'echet-Urysohn} property.

\item A space $X$ is said to be \emph{well-based} at $x$ if $x$ has a neighbourhood base well-ordered by $\supseteq$.
Such a neighbourhood base is said to be \emph{well-ordered neighbourhood base}.

A subspace of $X$ that contains $x$ and is well-based at $x$ is called a \emph{spoke at $x$}.
We will denote the set of (closed) spokes at $x$ by $\spokes{x}{X}$ ($\closedSpokes{x}{X}$).

\item For a point $x$ in a space $X$, we denote its neighbourhood filter by $\nhoodFilterSup{x}{X}$, or $\nhoodFilter{x}$ when the space is unambiguous.
We define its \emph{neighbourhood core} to be the intersection of all neighbourhoods of $x$.
This will be denoted by $\nhoodCoreSup{x}{X}$, or $\nhoodCore{x}$ again if $X$ is unambiguous.
Note that in a $\Tone$-space, $\nhoodCore{x} = \set{x}$.
\end{itemize}
\end{definition}

The reason to introduce these strictly convergent sequences is because they allow us to construct spokes:

\begin{lemma}\cite[Claim in Theorem 4.1, pg.~16]{Leek_internal_2014}
Let $X$ be a topological space, $x \in X$ be given and let $f : \lambda \to X$ be an injective transfinite sequence that converges strictly to $x$.
Then $\nhoodCore{x} \cup \range{f}$ is a spoke at $x$.
\end{lemma}

We refer to spokes of this form as \emph{basic spokes}.
In \cite{Leek_internal_2014}, we used spoke systems consisting of basic spokes in our proof characterising radiality.

\begin{definition}[Spoke system]
Let $X$ be a topological space, $x \in X$ be given and let $\mathcal{S}$ be a collection of spokes at $x$.
Then we say that $\mathcal{S}$ is a \emph{spoke system} at $x$ if
\begin{equation*}
\bigset{\bigcup_{S \in \mathcal{S}} U_S : \forall S \in \mathcal{S}, U_S \in \nhoodFilterSup{x}{S}}
\end{equation*}
is a neighbourhood base for $x$ with respect to $X$.
Note that this collection will always form a network at $x$.
\end{definition}

\begin{definition}[Almost-independent]
Let $X$ be a topological space, $x \in X$ be given and let $\mathcal{S}$ be a collection of spokes at $x$ such that $x \notin \closure{(S \cap T) \backslash \nhoodCore{x}}$ for all distinct $S, T \in \mathcal{S}$.
Then we say that $\mathcal{S}$ is \emph{almost-independent}.
\end{definition}

\begin{theorem} \cite[Theorem 4.1, pg.~16]{Leek_internal_2014} \label{Spoke characterisation}
Let $X$ be a topological space, $x \in X$ be given.
Then the following are equivalent:
\begin{enumerate}
\item $X$ is radial at $x$.

\item $x$ has a spoke system.

\item $x$ has an almost-independent, basic spoke system.
\end{enumerate}
\end{theorem}

If we assume some extra separation axioms, we can thicken our spokes.
This process will be useful when investigating radiality in compact spaces, and in particular compactifications.

\begin{definition}[Regular]
A point in a topological space is \emph{regular} if it has a neighbourhood base consisting of closed sets.
\end{definition}

\begin{lemma}
Let $x$ be a regular point in a topological space $X$ and let $S$ be a spoke of $x$.
Then $\closure{S}$ is also a spoke of $x$.
\begin{proof}
Choose a well-ordered neighbourhood base $\transSequence{B_\alpha}{\alpha}{\lambda}$ of $x$ with respect to $S$ and define for all $\alpha < \lambda, C_\alpha \df \closure{S} \backslash \closure{S \backslash B_\alpha} \in \nhoodFilterSup{x}{\closure{S}}$.
Note that for all $\alpha < \beta < \lambda, C_\beta \subseteq C_\alpha$.
Let $D$ be a closed neighbourhood of $x$ with respect to $X$, so there exists an $\alpha < \lambda$ such that $B_\alpha \subseteq D$ and hence $S \backslash D \subseteq S \backslash B_\alpha$.
Then
\begin{equation*}
S \subseteq D \cup \closure{S \backslash B_\alpha} \Rightarrow \closure{S} \subseteq D \cup \closure{S \backslash B_\alpha} \Rightarrow C_\alpha \subseteq D \cap \closure{S}.
\end{equation*}
Therefore $\transSequence{C_\alpha}{\alpha}{\lambda}$ is a well-ordered neighbourhood base of $x$ with respect to $\closure{S}$.
As $\nhoodCore{x} \subseteq S \subseteq \closure{S}$, it follows that $\closure{S}$ is a spoke of $x$.
\end{proof}
\end{lemma}

\begin{corollary} \label{Compact Hausdorff radial characterisation}
Let $X$ be a regular space and let $x \in X$ be given.
Then $X$ is radial at $x$ if and only if $x$ has a closed spoke system.
\begin{proof}
By Theorem \ref{Spoke characterisation}, it suffices to assume $X$ has a closed spoke system $\mathcal{S}$ and define $\mathcal{T} \df \set{\closure{S} : S \in \mathcal{S}}$, which by the previous lemma is a collection of spokes of $x$.
For each $T \in \mathcal{T}$, choose a $U_T \in \mathcal{N}_x$ and define for all $S \in \mathcal{S}, V_S \df U_{\closure{S}}$.
Then
\begin{equation*}
\bigcup_{T \in \mathcal{T}} (U_T \cap T) \supseteq \bigcup_{S \in \mathcal{S}} (V_{\closure{S}} \cap S) \in \mathcal{N}_x^X.
\end{equation*}
Therefore $\mathcal{T}$ is a closed spoke system of $x$.
\end{proof}
\end{corollary}

However, we don't necessarily have a spoke system that is both closed and almost-independent, even for compact Hausdorff spaces.
We need to introduce some more notation: we will denote the one-point compactification of a space $X$ by $\alpha X$, with it's point-at-infinity denoted by $\star$.
Also, let $\compactSubsets{X}$ denote the set of compact subsets of a topological space $X$.

\begin{theorem}
There exists a compact Hausdorff space $X$ and a radial point $x \in X$ with no closed, almost-independent spoke system.
\begin{proof}
Define $X \df \alpha(\omega_1 \times \omega_2)$ and note that for all $K \in \compactSubsets{\omega_1 \times \omega_2}, \pi_{\omega_1}[K], \pi_{\omega_2}[K]$ are bounded in $\omega_1, \omega_2$ respectively and hence $K \subseteq \alpha \times \beta$ for some $\alpha < \omega_1$ and $\beta < \omega_2$.
In particular, every $\sigma$-compact subset of $\omega_1 \times \omega_2$ has compact closure; i.e., $\star$ is a p-point\footnote{A point $x$ in a topological space is a \emph{p-point} if countable intersections of neighbourhoods of $x$ are again a neighbourhood.}.

Let $A \subseteq \omega_1 \times \omega_2$ be given such that $\star \in \closure{A}$.
Then $\closure{A}^{\omega_1 \times \omega_2}$ is not compact, so there exists an $i = 1, 2$ such that $\pi_{\omega_i}[A]$ is unbounded in $\omega_i$.
Then for all $\alpha < \omega_i$, there exists an $a_\alpha \in A$ such that $\pi_{\omega_i}(a_\alpha) > \alpha$.
Let $K \in \compactSubsets{\omega_1 \times \omega_2}$ be given, so there exists an $\alpha_1 < \omega_1$ and $\alpha_2 < \omega_2$ such that $K \subseteq \alpha_1 \times \alpha_2$.
Then $a_\beta \notin K$ for all $\beta \in [\alpha_i, \omega_i)$, so $\transSequence{a_\beta}{\beta}{\omega_i} \to \star$.
Therefore $\star$ is radial in $X$.

Now suppose there exists a closed, almost-independent spoke system $\mathcal{S}$ for $\star$ and define $\Lambda \df \set{\lambda < \omega_2 : \cofinality{\lambda} = \omega_1}$.
We claim that for all $\lambda \in \Lambda$, there exists an $\alpha_\lambda < \omega_1$, a $\beta_\lambda < \lambda$ and an $S_\lambda \in \mathcal{S}$ such that $[\alpha_\lambda, \omega_1) \times [\beta_\lambda, \lambda) \subseteq S_\lambda$.
Before proving this claim, we will show how it will allow us to derive a contradiction.

Suppose that $\set{S_\lambda : \lambda \in \Lambda}$ is uncountable and pick $f : \omega_1 \to \Lambda$ strictly increasing such that for all distinct $\alpha, \beta < \omega_1, S_{f(\alpha)} \ne S_{f(\beta)}$.
Define $\lambda \df \sup(\range{f}) \in \Lambda$.
Then $f$ is cofinal in $\lambda$, so there exists a $\gamma < \omega_1$ such that $f(\gamma) > \beta_\lambda$.
Thus for all $\delta \in [\gamma, \omega_1)$:
\begin{equation*}
[\max(\alpha_{f(\delta)}, \alpha_\lambda), \omega_1) \times \set{\max(\beta_{f(\delta)}, \beta_\lambda)} \subseteq S_{f(\delta)} \cap S_\lambda
\end{equation*}
Hence $\star \in \closure{(S_{f(\delta)} \cap S_\lambda) \backslash \set{\star}}$.
Since $\mathcal{S}$ is almost-independent, it follows that $S_{f(\delta)} = S_\lambda$ and in particular $S_{f(\gamma)} = S_{f(\gamma + 1)}$, which is a contradiction.
Therefore $\set{S_\lambda : \lambda \in \Lambda}$ is countable and so there exists an $L \subseteq \Lambda$ of cardinality $\aleph_2$ such that $S_\lambda = S_\mu$ for all $\lambda, \mu \in L$.
As each $S_\lambda$ contains a non-trivial $\omega_1$-sequence converging to $\star$, it follows that $\chi(\star, S_\lambda) = \aleph_1$.
However, $\family{(\alpha_\lambda, \beta_\lambda)}{\lambda}{L}$ is an $\omega_2$-sequence in $S_{\min(L)}$ that converges to $\star$, which is a contradiction.
Therefore $\star$ doesn't have a closed, almost-independent spoke system.

We will now prove our claim.
Let $f : \omega_1 \to \omega_1 \times \lambda$ be given such that $f \to \star$.
Suppose for all $S \in \mathcal{S}$, there exists a $U_S \in \nhoodFilterSup{\star}{S}$ such that $\range{f} \cap U_S = \emptyset$.
Then $U \df \bigcup_{S \in \mathcal{S}} U_S \in \nhoodFilterSup{\star}{X}$ and $\range{f} \cap U = \emptyset$, which is a contradiction.
Thus there exists an $S \in \mathcal{S}$ such that $\star \in \closure{\range{f} \cap S}$.
Since $\star$ is a p-point, it follows that $\range{f} \cap S$ is uncountable.
Now let $h : \omega_1 \to \lambda$ be cofinal and strictly increasing and continuous.
Define for all $\alpha < \omega_1, f(\alpha) \df (\alpha, h(\alpha))$.
Then $f \to \star$, so by the work above there exists an $S_\lambda \in \mathcal{S}$ such that $\range{f} \cap S_\lambda$ is uncountable.
Define $A \df \pi_{\omega_1}[\range{f} \cap S_\lambda]$.

Suppose for all $\alpha < \omega_1$, there exists an $x_\alpha \in ([\alpha, \omega_1) \times [h(\alpha), \lambda)) \backslash S_\lambda$.
Then $\transSequence{x_\alpha}{\alpha}{\omega_1} \to \star$, so again there exists a $T \in \mathcal{S}$ such that $B \df \pi_{\omega_1}[\set{x_\alpha : \alpha < \omega_1} \cap T]$ is uncountable.
Since $\set{x_\alpha : \alpha < \omega_1} \cap S_\lambda = \emptyset$, it follows that $T$ is distinct from $S_\lambda$, so $(S_\lambda \cap T) \backslash \set{\star}$ has compact closure in $\omega_1 \times \omega_2$.
In particular, its projection onto $\omega_1$ is bounded.

Let $\beta < \omega_1$ be given.
As $\lambda$ has uncountable cofinality, $A' \cap B'$ is a club.
Let $\gamma \in [\beta, \omega_1) \cap A' \cap B'$ be given, so there exist strictly increasing sequences $\transSequence{\delta_n}{n}{\omega} \subseteq A, |\transSequence{\epsilon_n}{n}{\omega} \subseteq B$ with supremum $\gamma$.
Then by continuity of $h, (\gamma, h(\gamma)) \in \closure{\set{(\delta_n, h(\delta_n)) : n < \omega}} \subseteq \closure{\range{f} \cap S_\lambda} \subseteq S_\lambda$.
Moreover, for each $n < \omega$, there exists an $\alpha_n < \omega_1$ such that $x_{\alpha_n} \in T$ and $\pi_{\omega_1}(x_{\alpha_n}) = \epsilon_n$.
Therefore, since $\lambda$ is sequentially compact, by virtue of being an ordinal with uncountable cofinality, there exists a subsequence of $\transSequence{\pi_{\omega_2}(x_{\alpha_n})}{n}{\omega_1}$ that converges to some ordinal $\theta < \lambda$ and so $(\gamma, \theta) \in \closure{\set{x_\alpha : \alpha \in B} \cap T} \subseteq T$.
Hence $\gamma \in \pi_{\omega_1}[(S_\lambda \cap T) \backslash \set{\star}]$.
But this then shows that $\pi_{\omega_1}[(S_\lambda \cap T) \backslash \set{\star}]$ is unbounded, which is a contradiction.
Thus there exists an $\alpha_\lambda < \omega_1$ such that $[\alpha_\lambda, \omega_1) \times [h(\alpha_\lambda), \lambda) \subseteq S_\lambda$.
By defining $\beta_\lambda \df h(\alpha_\lambda)$, we conclude the proof of our claim and the theorem.
\end{proof}
\end{theorem}

To finish this section, we will exhibit an alternative characterisation of radiality by ordering our spokes by local containment.
This has the added advantage of characterising the \emph{subspaces} which are radial at a specified point.

\begin{definition}[Locally contained]
Let $X$ be a topological space, $x \in X, A, B \subseteq X$ be given.
Then we say \emph{locally at $x, A$ is contained in $B$}, written $A \locallyContainedSup{x}{X} B$, if there exists a $U \in \nhoodFilter{x}$ such that $A \cap U \subseteq B$, or equivalently, $x \notin \closure{A \backslash B}$.
If the ambient space $X$ is unambiguous, we will drop the superscript in $\locallyContainedSup{x}{X}$.
\end{definition}

We will endow $\spokes{x}{X}$ and $\closedSpokes{x}{X}$ with this ordering and consider cofinal subsets of this quasi-ordered set.

\begin{lemma}
Let $X$ be a topological space and let $Y \subseteq X, x \in Y, S \in \spokes{x}{Y}$ be given.
Then $S \cup \nhoodCoreSup{x}{X} \in \spokes{x}{X}$.
\begin{proof}
Let $\set{B_\alpha : \alpha < \lambda}$ be a well-ordered neighbourhood base for $x$ with respect to $Y$, where for all $\alpha, \beta < \lambda$ with $\alpha < \beta, B_\beta \subseteq B_\alpha$.
Let $U \in \nhoodFilterSup{x}{X}$ be given, so there exists an $\alpha < \lambda$ such that $B_\alpha \subseteq U \cap S$.
Then $B_\alpha \cup \nhoodCoreSup{x}{X} \subseteq U \cap (S \cup \nhoodCoreSup{x}{X})$.
Moreover, for all $\alpha < \lambda$, there exists a $V \in \nhoodFilter{x}{X}$ such that $B_\alpha = V \cap S$ and so $B_\alpha \cup \nhoodCoreSup{x}{X} = V \cap (S \cup \nhoodCoreSup{x}{X}) \in \nhoodFilterSup{x}{S \cup \nhoodCoreSup{x}{X}}$.
Therefore $\set{B_\alpha \cup \nhoodCoreSup{x}{X} : \alpha < \lambda}$ is a well-ordered neighbourhood base for $x$ with respect to $S \cup \nhoodCoreSup{x}{X}$.
Hence $S \cup \nhoodCoreSup{x}{X} \in \spokes{x}{X}$.
\end{proof}
\end{lemma}

\begin{theorem} \label{Radial subspace locally}
Let $X$ be a topological space, $Y \subseteq X, x \in Y$ be given.
Then the following are equivalent:
\begin{enumerate}
\item $Y$ is radial at $x$.

\item For all $\mathcal{C} \subseteq \spokes{x}{X}$ cofinal, $Y \locallyContainedSup{x}{X} \bigcup \mathcal{C}$.
\end{enumerate}
Moreover, if $X$ is compact and Hausdorff, then the two conditions above are equivalent to:
\begin{enumerate}
\setcounter{enumi}{2}
\item For all $\mathcal{C} \subseteq \closedSpokes{x}{X}$ cofinal, $Y \locallyContainedSup{x}{X} \bigcup \mathcal{C}$.
\end{enumerate}
\begin{proof}
Suppose $Y$ is radial at $x$, so there exists a spoke system $\mathcal{S}$ for $x$ with respect to $Y$.
Let $\mathcal{C} \subseteq \spokes{x}{X}$ be cofinal, so by the previous lemma for all $S \in \mathcal{S}$, there exists a $C_S \in \mathcal{C}$ such that $S \cup \nhoodCoreSup{x}{X} \locallyContainedSup{x}{X} C_S$ and thus there exists a $U_S \in \nhoodFilterSup{x}{X}$ such that $(S \cap U_S) \cup \nhoodCoreSup{x}{X} = (S \cup \nhoodCoreSup{x}{X}) \cap U_S \subseteq C_S$.
Define:
\begin{equation*}
U \df \bigcup_{S \in \mathcal{S}} (S \cap U_S) \in \nhoodFilterSup{x}{Y}
\end{equation*}
Then there exists a $V \in \nhoodFilterSup{x}{X}$ such that $U = V \cap Y$.
As $U \subseteq \bigcup \mathcal{C}$, it follows that $Y \locallyContainedSup{x}{X} \bigcup \mathcal{C}$.
Thus (1) implies (2).
Moreover, if $X$ is compact and Hausdorff, then we can take $\mathcal{S}$ to consist of closed spokes by Corollary \ref{Compact Hausdorff radial characterisation} and $\mathcal{C} \subseteq \closedSpokes{x}{X}$.
Therefore (1) implies (3) too.

Now suppose that for all $\mathcal{C} \subseteq \spokes{x}{X}$ cofinal, $Y \locallyContainedSup{x}{X} \bigcup \mathcal{C}$.
We will show that $\set{S \cap Y : S \in \spokes{x}{X}}$ is a spoke system for $x$ with respect to $Y$.
For all $S \in \spokes{x}{X}$, let $U_S \in \nhoodFilterSup{x}{X}$ be given, so $U_S \cap Y \in \nhoodFilterSup{x}{Y}$.
Note that $\set{S \cap U_S : S \in \spokes{x}{X}}$ is vacuously cofinal in $\spokes{x}{X}$, so there exists a $V \in \nhoodFilterSup{x}{X}$ such that $V \cap Y \subseteq \bigcup_{S \in \mathcal{S}} (S \cap U_S)$.
Then $V \cap Y \subseteq \bigcup_{S \in \mathcal{S}} ((S \cap Y) \cap (U_S \cap Y))$, so the latter is a neighbourhood of $x$ with respect to $Y$.
Therefore $\set{S \cap Y : S \in \spokes{x}{X}}$ is a spoke system for $x$ with respect to $Y$ and thus $Y$ is radial at $x$ by Theorem \ref{Spoke characterisation}.
Hence (2) implies (1).
Finally, note that if $X$ is compact and Hausdorff then by replacing $\spokes{x}{X}$ with $\closedSpokes{x}{X}$, we see that (3) implies (1), concluding our proof.
\end{proof}
\end{theorem}

\begin{corollary} \label{Cofinal spoke collections}
Let $X$ be a topological space, $x \in X$ be given.
Then $X$ is radial at $x$ if and only if for all cofinal collections of spokes $\mathcal{C}, \bigcup \mathcal{C}$ is a neighbourhood of $x$.
\begin{proof}
By the previous theorem, $X$ is radial at $x$ if and only if $X \locallyContained{x} \bigcup \mathcal{C}$ for all cofinal collections $\mathcal{C} \subseteq \spokes{x}{X}$, which is equivalent $\bigcup \mathcal{C} \in \nhoodFilter{x}$.
\end{proof}
\end{corollary}

\section{One-point compactifications}

For the rest of this article, unless otherwise stated, we will assume that $X$ is a locally compact, non-compact Hausdorff space.

In this section, we will use our spoke characterisations to characterise being radial at $\star$ in $\alpha X$.
Spoke systems and cofinal collections of spokes allow us to reflect these properties from compactifications down to the structure of the compact subsets of $X$.
We will first show that it suffices to consider the points in the remainder.

\begin{lemma} \label{Radial compactification reduction}
Let $X$ be a topological space, $U \subseteq X$ be open.
If $U$ is radial and $X$ is radial at every point outside $U$ then $X$ is radial.
\begin{proof}
Let $u \in U, A \subseteq X$ be given such that $u \in \closure{A}$.
Then for each $V \subseteq X$ open, if $u \in V$ then $V \cap A \ne \emptyset$.
In particular, for each $V \subseteq U$ open, $V \cap A \ne \emptyset$ and so $u \in \closure{A \cap U}^U$.
Thus there exists a transfinite sequence contained in $A \cap U$ that converges to $u$ and therefore $X$ is radial at $u$.
\end{proof}
\end{lemma}

\begin{definition}[Spoke at infinity]
Let $S \subseteq X$ be given such that $S \cup \set{\star}$ is a spoke at $\star$ in $\alpha X$.
Then we say that $S$ is a \emph{spoke at infinity} of $X$.
We will denote the set of (closed) spokes at infinity by $\spokesAtInfinity{X}$ ($\closedSpokesAtInfinity{X}$).
\end{definition}

\begin{lemma} \label{Spoke infinity characterisation}
Let $S \subseteq X$ be closed.
Then $S$ is a spoke at infinity if and only if there exists a cofinal chain in $(\compactSubsets{S}, \subseteq)$.
\begin{proof}
Assume $S$ is a spoke at infinity, so there exists a well-ordered neighbourhood base $\mathcal{B}$ of $\star$ in $S \cup \set{\star}$.
By taking interiors, we can assume that $\mathcal{B}$ consists of open sets.
Then $\set{S \backslash B : B \in \mathcal{B}}$ is a chain in $\compactSubsets{S}$.
Moreover, for all $K \in \compactSubsets{S}$, there exists a $B \in \mathcal{B}$ such that $B \subseteq \alpha X \backslash K$ and so $K \subseteq S \backslash B$.
Therefore $\set{S \backslash B : B \in \mathcal{B}}$ is cofinal in $(\compactSubsets{S}, \subseteq)$.

Now assume that there exists a cofinal chain in $(\compactSubsets{S}, \subseteq)$, so by considering its cofinality, there exists an increasing, cofinal, transfinite sequence $\transSequence{K_\alpha}{\alpha}{\lambda} \subseteq \compactSubsets{S}$.
Then $(S \cup \set{\star}) \backslash K_\alpha$ is a neighbourhood of $\star$ in $S$ for each $\alpha < \lambda$.
Let $V \subseteq \alpha X$ be open with $\star \in V$, so $X \backslash V$ is compact and hence $S \backslash V$ has compact closure in $S$.
Thus there exists an $\alpha < \lambda$ such that $S \backslash V \subseteq K_\alpha$ and so $(S \cup \set{\star}) \backslash K_\alpha \subseteq S \cap V$.
Therefore $\transSequence{(S \cup \set{\star}) \backslash K_\alpha}{\alpha}{\lambda}$ is a well-ordered neighbourhood base for $\star$ in $S \cup \set{\star}$ and hence $S$ is a spoke at infinity.
\end{proof}
\end{lemma}

The following theorem demonstrates an internal characterisation for radiality at infinity, using the spoke system criterion.

\begin{theorem} \label{Alexandroff compactification radial characterisation}
$\alpha X$ is radial at $\star$ if and only if there exists a collection $\mathcal{S} \subseteq \closedSpokesAtInfinity{X}$ such that for all $C \in \prod_{S \in \mathcal{S}} \compactSubsets{S}, \bigcup_{S \in \mathcal{S}} (S \backslash C(S))$ has co-compact\footnote{A subset is \emph{co-compact} if its complement is compact.} interior.
\begin{proof}
Let $\mathcal{S} \subseteq \closedSpokes{\star}{\alpha X}$ be given.
Then $\mathcal{S}$ is a spoke system at $\star$ if and only if for all $C \in \prod_{S \in \mathcal{S}} \nhoodFilterSup{\star}{S}, \linebreak
\bigcup_{S \in \mathcal{S}} C(S) \in \nhoodFilterSup{x}{\alpha X}$.
Since $\alpha X$ is compact, this is equivalent to $\bigcup_{S \in \mathcal{S}} ((S \backslash \set{\star}) \backslash C(S))$ having co-compact interior in $X$ for all $C \in \prod_{S \in \mathcal{S}} \compactSubsets{S \backslash \set{\star}}$.
Thus by Corollary \ref{Compact Hausdorff radial characterisation}, the proof is complete.
\end{proof}
\end{theorem}

We also have the following characterisation in terms of spokes at infinity, purely from the radiality property itself.

\begin{theorem} \label{Radial trace spoke at infinity}
$\alpha X$ is radial at $\star$ if and only if for all $Y \subseteq X$ with non-compact closure in $X$, there exists a non-compact $Z \in \closedSpokesAtInfinity{X}$ such that $K_\alpha = \closure{K_\alpha \cap Y}$ for all $\alpha < \lambda$, where $\transSequence{K_\alpha}{\alpha}{\lambda}$ is a cofinal chain in $(\compactSubsets{Z}, \subseteq)$.
\begin{proof}
Suppose that $\alpha X$ is radial at $\star$ and let $Y \subseteq X$ have non-compact closure in $X$, so $\star \in \closure{Y}$.
Then by radiality there exists an injective transfinite sequence $f : \lambda \to Y$ that converges strictly to $\star$ (see \cite[Lemma 2.2, pg.~12]{Leek_internal_2014}), so $S(f) \in \closedSpokesAtInfinity{X}$ and $\transSequence{\closure{f[\beta]}}{\beta}{\lambda}$ is a cofinal chain in $(\compactSubsets{S(f)}, \subseteq)$.
Let $\beta < \lambda$ be given.
Then:
\begin{equation*}
f[\beta] \subseteq \closure{f[\beta]} \cap Y \subseteq \closure{\closure{f[\beta]} \cap Y} \subseteq \closure{f[\beta]}
\end{equation*}
Hence $\closure{f[\beta]} = \closure{\closure{f[\beta]} \cap Y}$.
Also, $S(f)$ is non-compact, since $\star \in \closure{S(f)}^{\alpha X}$.

Now suppose the converse holds and let $A \subseteq X$ be given such that $\star \in \closure{A}$, so $A$ has non-compact closure.
Then there exists a non-compact $Y \in \closedSpokesAtInfinity{X}$, with $\transSequence{K_\beta}{\beta}{\lambda}$ a strictly increasing, cofinal chain in $\compactSubsets{Y}$, such that $K_\beta = \closure{K_\beta \cap A}$ for all $\beta < \lambda$.
Since $Y$ is non-compact, $\lambda$ must be a limit ordinal.
Let $\beta < \lambda$ be given, so $\closure{K_\beta \cap A} = K_\beta \subsetneqq K_{\beta + 1} = \closure{K_{\beta + 1} \cap A}$ and hence there exists an $x_\beta \in (K_{\beta + 1} \backslash K_\beta) \cap A$.
Now since $\transSequence{\set{\star} \cup (Y \backslash K_\beta)}{\beta}{\lambda}$ is a neighbourhood base for $\star$ with respect to $Y \cup \set{\star}$ (by the proof of Lemma \ref{Spoke infinity characterisation}), it follows that $\transSequence{x_\beta}{\beta}{\lambda}$ converges to $\star$ and is contained in $A$.
Therefore $\alpha X$ is radial at $\star$.
\end{proof}
\end{theorem}

\begin{corollary}
Suppose $\alpha X$ is radial at $\star$.
Then for all $A \subseteq X$ closed and non-compact, there exists a non-compact $S \in \closedSpokesAtInfinity{X}$ contained in $A$.
\begin{proof}
By picking $S$ and $\transSequence{K_\alpha}{\alpha}{\lambda}$ from the previous theorem, it follows that $K_\alpha = \closure{K_\alpha \cap A} \subseteq A$ for all $\alpha < \lambda$ and so $S = \bigcup_{\alpha < \lambda} K_\alpha \subseteq A$.
\end{proof}
\end{corollary}

Of course, the preceding corollary is not surprising when our spokes at infinity are $\sigma$-compact, for we can then take an $\omega$-sequence converging to $\star$.
However, this is more an artefact of $\Tone$ implying finite subsets are closed.
If $\star$ is a p-point, then the spokes will contain closures of countably-infinite subsets, which could potentially be large.
Unfortunately, even though this corollary is a more natural condition, it is not equivalent to radiality at $\star$:

\begin{theorem}
There exists a non-compact, locally compact Hausdorff space such that for all $A \subseteq X$ closed and non-compact, there exists a non-compact $S \in \closedSpokesAtInfinity{X}$ with $S \subseteq A$, yet $\alpha X$ is not radial at $\star$.
\begin{proof}
Define the \emph{deleted Tychonoff plank} to be $X \df ((\omega + 1) \times (\omega_1 + 1)) \backslash \set{(\omega, \omega_1)}$ and observe that $\alpha X \cong (\omega + 1) \times (\omega_1 + 1)$, so $\alpha X$ is not radial at $\star = (\omega, \omega_1)$ (as noted in \cite[pg.~12-13]{Leek_internal_2014}). Let $A \subseteq X$ be closed and non-compact and suppose $\pi_\omega[A \cap (\omega \times \set{\omega_1})]$ and $\pi_{\omega_1}[A \cap (\set{\omega} \times \omega_1)]$ are bounded in $\omega$ and $\omega_1$ respectively.
Then there exists an $n \in \omega$ and $\alpha < \omega_1$ such that $A \subseteq ((\omega + 1) \times (\omega_1 + 1)) \backslash ([n, \omega] \times [\alpha, \omega_1])$.
Thus $A \subseteq (n \times (\omega_1 + 1)) \cup ((\omega + 1) \times (\alpha + 1))$, which is compact and hence a contradiction.
Therefore either $\pi_\omega[A \cap (\omega \times \set{\omega_1})]$ is unbounded in $\omega$ or $\pi_{\omega_1}[A \cap (\set{\omega} \times \omega_1)]$ is unbounded in $\omega_1$.
As $\set{\omega} \times \omega_1$ and $\omega \times \set{\omega_1}$ are easily seen to be spokes at infinity, it follows that $A \cap (\omega \times \set{\omega_1}) \in \closedSpokesAtInfinity{X}, A \cap (\set{\omega} \times \omega_1) \in \closedSpokesAtInfinity{X}$ and one of these is non-compact.
This completes the proof.
\end{proof}
\end{theorem}

We will now present the third characterisation using cofinal spoke collections.
We first need to translate the ordering on spokes of $\star$ to spokes at infinity.

\begin{definition}
Let $X$ be a non-compact, locally compact Hausdorff space.
For all $S, T \in \spokesAtInfinity{X}$, we define $S \leq T$ if $S \backslash T$ has compact closure.
Observe that for $S, T \in \spokesAtInfinity{X}, S \cup \set{\star} \locallyContainedSup{\star}{\alpha X} T \cup \set{\star}$ if and only if $S \leq T$.
We will endow $\spokesAtInfinity{X}$ and $\closedSpokesAtInfinity{X}$ with this quasi-order.
\end{definition}

\begin{theorem}
$\alpha X$ is radial at $\star$ if and only if for all $\mathcal{C} \subseteq \closedSpokesAtInfinity{X}$ cofinal, $\bigcup \mathcal{C}$ has co-compact interior.
\begin{proof}
Assume $\alpha X$ is radial at $\star$ and let $\mathcal{C} \subseteq \closedSpokesAtInfinity{X}$ be cofinal.
Then for all $S \in \closedSpokesAtInfinity{X}$, there exists a $T_S \in \mathcal{C}$ such that $C_S \df \closure{S \backslash T_S}$ is compact.
As $\alpha X$ is radial at $\star$, it follows by Theorem \ref{Alexandroff compactification radial characterisation} that $\bigcup_{S \in \closedSpokesAtInfinity{X}} (S \backslash C_S)$ has co-compact interior.
Note that for all $S \in \closedSpokesAtInfinity{X}, S \backslash C_S \subseteq S \backslash (S \backslash T_S) \subseteq T_S$.
Thus $\bigcup_{S \in \closedSpokesAtInfinity{X}} (S \backslash C_S) \subseteq \bigcup \mathcal{C}$, so $\bigcup \mathcal{C}$ also has co-compact interior.

Now assume that for all $\mathcal{C} \subseteq \closedSpokesAtInfinity{X}$ cofinal, $\bigcup \mathcal{C}$ has co-compact interior.
Let $C \in \prod_{S \in \closedSpokesAtInfinity{X}} \compactSubsets{S}$ be given and define $\mathcal{C} \df \set{S \backslash C(S) : S \in \closedSpokesAtInfinity{X}}$.
Then $\mathcal{C}$ is cofinal in $\closedSpokesAtInfinity{X}$ and hence $\bigcup \mathcal{C}$ has co-compact interior.
Therefore by Theorem \ref{Alexandroff compactification radial characterisation} again, $\alpha X$ is radial at $\star$.
\end{proof}
\end{theorem}

We will now analyse two spaces, which are known to not be radial at infinity, and proving this fact using these theorems.

\subsection{Deleted Tychonoff plank}

Let $P$ denote the deleted Tychonoff plank and as before we may take $\alpha P = (\omega + 1) \times (\omega_1 + 1)$ and $\star = (\omega, \omega_1)$.
Define $S_0 \df \omega \times \set{\omega_1}, S_1 \df \set{\omega} \times \omega_1$.
Note that in an ordinal space, every compact subset is bounded and $S_0 = \bigcup_{n \in \omega} (n \times \set{\omega_1}), S_1 = \bigcup_{\alpha < \omega_1} (\set{\omega} \times (\alpha + 1))$, so $S_0, S_1 \in \closedSpokesAtInfinity{X}$.
We will show that $\set{S_0, S_1}$ is a cofinal collection of closed spokes at infinity.

Let $S \in \closedSpokesAtInfinity{P}$ be non-compact and let $\transSequence{K_\alpha}{\alpha}{\lambda} \subseteq \compactSubsets{S}$ be a cofinal chain with $\lambda$ infinite.
Without loss of generality, assume $\lambda$ is regular and for all $\alpha < \lambda, K_\alpha \subsetneqq K_{\alpha + 1}$.
Since $|P| = \aleph_1$, either $\lambda = \omega$ or $\lambda = \omega_1$.
\begin{description}
\item[Case 1:] Suppose $\lambda = \omega$ and consider $S \cap S_1$.
Then $(\compactSubsets{S \cap S_1}, \subseteq)$ has cofinal chains of lengths $\omega$ and $\omega_1$, so $S \cap S_1$ must be compact and hence there exists a $\beta < \omega_1$ such that $S \cap S_1 \subseteq \set{\omega} \times \beta$.
Assume $S \not \leq S_0$, so for all $n < \omega$, there exists an
\begin{equation*}
x_n \in S \backslash ((n \times (\omega_1 + 1)) \cup ((\omega + 1) \times (\beta + 1)) \cup S_0) =  S \cap (((\omega + 1) \backslash n) \times (\omega_1 \backslash \beta)).
\end{equation*}
Since $\omega_1$ is sequentially compact, there exists a strictly increasing sequence $\transSequence{r_n}{n}{\omega}$ in $\omega$ and $\gamma \in \omega_1 \backslash \beta$ such that $\transSequence{\pi_{\omega_1}(x_{r_n})}{n}{\omega} \to \gamma$ and hence $\transSequence{x_{r_n}}{n}{\omega} \to (\omega, \gamma)$.
However, since $S$ is closed, $(\omega, \gamma) \in S \cap S_1$, which is a contradiction.
Therefore $S \leq S_0$.

\item[Case 2:] Suppose $\lambda = \omega_1$ and consider $S \cap S_0$.
Again, $(\compactSubsets{S \cap S_0}, \subseteq)$ has cofinal chains of lengths $\omega$ and $\omega_1$, so $S \cap S_0$ must be compact and hence there exists an $n < \omega$ such that $S \cap S_0 \subseteq n \times \set{\omega_1}$.
Since $S$ is closed, for all $m \geq n$ there exists a $\beta_m < \omega_1$ such that $S \cap (\set{m} \times ((\omega + 1) \backslash \beta_m)) = \emptyset$.
Define $\beta \df \sup(\set{\beta_m : m \geq n}) < \omega_1$, so $S \cap ((\omega \backslash n) \times ((\omega_1 + 1) \backslash \beta)) = \emptyset$.
Then:
\begin{equation*}
S \backslash ((n \times (\omega_1 + 1)) \cup ((\omega + 1) \times (\beta + 1)) \cup S_1) = (S \backslash S_1) \cap (((\omega + 1) \backslash n) \times ((\omega_1 + 1) \backslash (\beta + 1))) = \emptyset
\end{equation*}
Therefore $S \leq S_1$.
\end{description}
Vacuously, every compact spoke at infinity is bounded above by $S_0$, so it follows that $\set{S_0, S_1}$ is a cofinal collection of paths to infinity.
However, $S_0 \cup S_1$ has empty, and hence non-co-compact, interior in $P$, so $\alpha P$ is not radial at $x$.

We also obtain a local result from Theorem \ref{Radial subspace locally}: since $\set{S_0, S_1}$ is a cofinal collection of closed spokes at infinity, any subspace of $(\omega + 1) \times (\omega_1 + 1)$ that is radial at $(\omega, \omega_1)$ must be locally contained at $(\omega, \omega_1)$ in $S_0 \cup S_1 \cup \set{(\omega, \omega_1)}$; indeed, as $S_0 \cup S_1 \cup \set{(\omega, \omega_1)}$ is a finite union of spokes at $(\omega, \omega_1)$, it is radial at $(\omega, \omega_1)$ and even a radial \emph{space}.

\subsection{Mr\'owka spaces}

Let $\mathcal{A}$ be a maximal, almost-disjoint (m.a.d.) family of subsets of $\omega$; that is, a maximal collection of infinite subsets of $\omega$ such that any two distinct elements intersect finitely.
We will define a topology on $\omega \cup \mathcal{A}$ as follows: let each $n \in \omega$ be isolated and for all $A \in \mathcal{A}$, let $\set{\set{A} \cup (A \backslash F) : F \subseteq \omega \text{ is finite}}$ be a neighbourhood base for $A$.
We denote this space by $\Psi$ and call it a \emph{Mr\'owka space}.
By \cite[Example 7.1, pg.~54-55]{Franklin_spaces_1967}, it is non-compact, locally compact and Hausdorff.
Moreover, it's one-point compactification is not radial at $\star$ since there is no (transfinite) sequence in $\omega$ converging to $\star$.
We will now show that the countably infinite subsets of $\mathcal{A}$ form a cofinal collection of closed spokes at infinity, witnessing this fact.

First note that $\mathcal{A}$ is closed and discrete in $\Psi$, so it easily follows that every countably infinite subset of $\mathcal{A}$ is a $\sigma$-compact spoke at infinity.
Let $S \in \closedSpokesAtInfinity{\Psi}$ be non-compact.
Then since no (transfinite) sequence in $\omega$ converges to $\star$, it follows that $S \cap \mathcal{A}$ must be a non-compact, closed spoke at infinity, since $\mathcal{A}$ is closed.
Since $\mathcal{A}$ is discrete, there are no infinite compact subsets of $\mathcal{A}$, so by Lemma \ref{Spoke infinity characterisation} $S \cap \mathcal{A}$ is countably infinite.
Thus $S = (S \cap \mathcal{A}) \cup (S \cap \omega)$ is countably infinite also.

Now suppose for all $\mathcal{F} \subseteq S \cap \mathcal{A}$ finite, $(S \cap \omega) \backslash (\bigcup \mathcal{F})$ is infinite and define
\begin{equation*}
\mathcal{B} \df \set{A \in S \cap \mathcal{A} : A \cap S \text{ is infinite}}.
\end{equation*}
Assume $\mathcal{B}$ is finite.
Then by maximality, there exists an $A \in \mathcal{A}$ such that $(A \cap S) \backslash (\bigcup \mathcal{B})$ is infinite and thus $A \in \closure{S} = S$, which is a contradiction.
Therefore $\mathcal{B}$ is infinite, so we can pick an enumeration $\mathcal{B} = \set{B_n : n < \omega}$.
Then for all $n < \omega$, there exists an $x_n \in (S \cap B_n) \backslash (\set{x_m : m < n} \cup \bigcup_{m < n} B_m)$ and so by maximality there is an $A \in \mathcal{A}$ such that $A \cap \set{x_n : n < \omega}$ is infinite.
Then we get a contradiction, since $A \in \closure{S} = S$ and $A \cap S$ is infinite, but $A \cap B_n$ is finite for all $n < \omega$.
Thus there exists a finite $\mathcal{F} \subseteq S \cap \mathcal{A}$ such that $C \df (S \cap \omega) \backslash (\bigcup \mathcal{F})$ is finite.
Then $S \cap \omega \subseteq C \cup \bigcup_{A \in \mathcal{F}} (\set{A} \cup A)$ and the latter is compact, so $S \leq S \cap \mathcal{A}$.
Therefore $[\mathcal{A}]^{\aleph_0} \df \set{\mathcal{A}' \subseteq \mathcal{A} : |\mathcal{A}'| = \aleph_0}$ is a cofinal collection of paths to infinity.
Note that $\bigcup [\mathcal{A}]^{\aleph_0} = \mathcal{A}$ has empty, and hence non-co-compact, interior in $\Psi$, so $\alpha \Psi$ is not radial at $\star$.

\section{Beyond the one-point compactification}

We will now investigate larger compactifications, assuming that $\alpha X$ is radial at $\star$.
We will start with finite and countable compactifications; in fact, we will demonstrate results for \emph{ordinal compactifications} - those which have remainder homeomorphic to some ordinal.
We also obtain conditions for the existence of sequential / pseudoradial compactifications.

Recall that we can obtain a one-point compactification of $X$ by identifying the remainder to a single point (see \cite[Theorems 3.5.12 \& 3.5.13, pg.~170]{engelking_general_1989}).
We will be implicitly using this identification from now on.

\subsection{Finite, countable and ordinal compactifications}

The following lemma demonstrates the usefulness of a space $X$ with the property that $\alpha X$ is radial at $\star$.

\begin{lemma} \label{Ordinal remainder subsequence lemma}
Let $\gamma X$ be an ordinal compactification of $X$; that is, $\gamma X \backslash X$ is homeomorphic to some ordinal and suppose $\alpha X$ is radial at $\star$.
Let $f : \lambda \to X$ be a transfinite sequence converges to $\star$ in $\alpha X$.
Then there is a subsequence\footnote{A \emph{subsequence} of a transfinite sequence $f : \lambda \to X$ is a transfinite sequence of the form $f \circ g : \mu \to X$, where $g : \mu \to \lambda$ is strictly increasing.} of $f$ that converges to some point in $\gamma X \backslash X$.
\begin{proof}
Assume not and pick an ordinal $\alpha$ such that $\gamma X \backslash X \cong \alpha + 1$ (since $\gamma X \backslash X$ is compact and non-empty).
We identify $\gamma X \backslash X$ with $\alpha + 1$.
By assumption, $\alpha$ must be non-zero.

Define $\lambda \df \domain{f}$ and suppose there exists an $m < \omega$ and a strictly decreasing sequence of ordinals $(\beta_n : n \leq m)$ in $\alpha + 1$ and $(U_n : n < m)$ a sequence of open subsets of $\gamma X$ such that:
\begin{itemize}
\item $\beta_0 = \alpha$,

\item $\beta_m > 0$,

\item $U_n \backslash X = (\beta_{n + 1}, \beta_n]$ for all $n < m$,

\item $D \df \lambda \backslash f^{-1}[\bigcup_{n < m} U_n]$ is unbounded.
\end{itemize}
Then since no subsequence of $f$ converges to $\beta_m$,  there exists an open subset $V \subseteq \gamma X$ such that $\beta_m \in V$ and $D \backslash f|_D^{-1}[V] = \lambda \backslash f^{-1}[V \cup \bigcup_{n < m} U_n]$ is unbounded.
Assume $[0, \beta_m) \subseteq V$.
Then $U \df V \cup \bigcup_{n < m} U_n$ is a neighbourhood of $\gamma X \backslash X$ and $\lambda \backslash f^{-1}[U]$ is unbounded, which is a contradiction since $f \to \star$ in $\alpha X$.
Thus there exists a $\beta_{m + 1} \in (0, \beta_m)$ such that $[\beta_{m + 1}, \beta_m] \subseteq V$ and furthermore there exists an open subset $W \subseteq \gamma X$ such that $(\beta_{m + 1}, \beta_m] = W \backslash X$.
Define $U_m \df V \cap W$ and note that $\lambda \backslash f^{-1}[\bigcup_{n \leq m} U_n]$ is unbounded and $(\beta_{m + 1}, \beta_m] = U_m \backslash X$.

Therefore by recursion, we find a descending sequence in $\alpha + 1$, which is a contradiction.
Hence there is a subsequence of $f$ that converges to some point in $\gamma X \backslash X$.
\end{proof}
\end{lemma}

Using the previous lemma, we can prove several compactification theorems.
Recall that a compactification is \emph{finite} / \emph{countable} if it has \emph{finite} / \emph{countable} remainder.

\begin{theorem}
Suppose $\alpha X$ is radial at $\star$ and let $\phi X$ be a finite compactification of $X$.
Then $\phi X$ is radial on $\phi X \backslash X$.
\begin{proof}
Let $A \subseteq \phi X, z \in \phi X \backslash X$ be given such that $z \in \closure{A}^{\phi X}$.
Since $\phi X$ is a finite compactification, there exists a closed neighbourhood $C \subseteq \phi X$ such that $\set{z} = C \backslash X$ and thus $z \in \closure{A \cap C}^{\phi X}$.
Then $\star \in \closure{A \cap C}^{\alpha X}$, so by radiality there exists a transfinite sequence $f$ contained in $A \cap C$ that converges $\star$.
By the previous lemma, $f$ has a subsequence that converges to some point in $\phi X \backslash X$.
This point must by $z$, and so $\phi X$ is radial on $\phi X \backslash X$.
\end{proof}
\end{theorem}

\begin{corollary} \label{Finite compactification}
Suppose $\alpha X$ is radial.
Then every finite compactification of $X$ is radial.
\begin{proof}
This follows from the previous theorem and Lemma \ref{Radial compactification reduction}.
\end{proof}
\end{corollary}

We now show that how we can obtain sequential / pseudoradial compactifications.
First, we recall their definitions:

\begin{definition}
Let $X$ be a topological space.
Then $X$ is \emph{pseudoradial} if for every non-closed subset $A \subseteq X$, there exists a transfinite sequence in $A$ that converges to a point outside $A$.
If we remove `transfinite' from this definiton, we obtain the definition of a \emph{sequential} space.
\end{definition}

\begin{lemma} \label{Pseudoradial compactification}
Suppose $\alpha X$ is radial / Fr\'echet-Urysohn at $\star$ and let $\gamma X$ be a compactification of $X$ with pseudoradial / sequential remainder such that every transfinite sequence / $\omega$-sequence in $X$ that converges to $\star$ in $\alpha X$ has a subsequence that converges to some point in $\gamma X \backslash X$.
Then $X$ is pseudoradial / sequential.
\begin{proof}
Let $A \subseteq \gamma X$ be non-closed and assume the first half of the conditions stated above.
Suppose $A \cap X$ is not closed in $X$.
Then since $X$ is pseudoradial, there is a transfinite sequence in $A \cap X$ that converges to a point in $X \backslash A$.
Now assume that $A \cap X$ is closed in $X$.
If $A \backslash X$ is not closed, then since $\gamma X \backslash X$ is pseudoradial, there exists a $z \in \closure{A \backslash X}^{\gamma X} \backslash A$ and a transfinite sequence contained in $A \backslash X$ that converges to $z$.

Finally, assume that $A \backslash X$ is closed.
Then since $A$ is not closed in $X$, there exists a $z \in \closure{A}^{\gamma X} \backslash A$ and thus $z \in \gamma X \backslash (A \cup X)$.
Since $A \backslash X$ is closed, there exists a closed neighbourhood $C \subseteq \gamma X$ of $z$ such that $C \cap (A \backslash X) = \emptyset$ and so $z \in \closure{A \cap C}^{\gamma X} = \closure{A \cap C \cap X}^{\gamma X}$.
Hence $\star \in \closure{A \cap C \cap X}^{\alpha X}$, so by radiality there exists a transfinite sequence $f$ in $A \cap C \cap X$ that converges to $\star$ in $\alpha X$ and thus by assumption there exists a $w \in \gamma X \backslash X$ and a subsequence $g$ of $f$ that converges to $w$.
As $C$ is closed in $\gamma X$, it follows that $w \in C$ and so $w \notin A$.
Therefore $\gamma X$ is pseudoradial.

For the second set of conditions, replace all occurrences of pseudoradial, radial and transfinite sequence with sequential, Fr\'echet-Urysohn and $\omega$-sequence respectively.
\end{proof}
\end{lemma}

\begin{theorem}
Suppose $X$ is pseudoradial and $\alpha X$ is radial at $\star$.
Then every ordinal compactification of $X$ is pseudoradial.
\begin{proof}
Note that every ordinal is (pseudo)radial, so by Lemmas \ref{Ordinal remainder subsequence lemma} and \ref{Pseudoradial compactification}, we have our result.
\end{proof}
\end{theorem}

\begin{theorem}
Suppose $X$ is pseudoradial / sequential and $\alpha X$ is radial / Fr\'echet-Urysohn at $\star$.
Then every compactification of $X$ with countable remainder is pseudoradial / sequential.
\begin{proof}
Note that every countable, compact Hausdorff space is homeomorphic to an ordinal (see \cite[pg.~351]{soukup_scattered_2004}), so is sequential.
Thus by Lemmas \ref{Ordinal remainder subsequence lemma} and \ref{Pseudoradial compactification}, we have our result.
\end{proof}
\end{theorem}

\subsection{Using small cardinals}

In this section, we will use \emph{small cardinals} to improve the results from the last section - these are uncountable cardinals bounded above by $\mathfrak{c} \df 2^{\aleph_0}$.
The small cardinals we are using are defined below.

\begin{definition}[Small cardinals] \hfill
\begin{itemize}
\item For $f, g : \omega \to \omega$, we say that $f$ is \emph{eventually bounded by $g$} if $\set{n \in \omega : f(n) > g(n)}$ is finite.
We denote this relation by $f \leq^* g$.
Observe that $\leq^*$ is a quasi-order on $\functions{\omega}{\omega}$.

\item The \emph{bounding number}, denoted by $\mathfrak{b}$, is the smallest cardinality of an unbounded subset of $(\functions{\omega}{\omega}, \leq^*)$.

\item For $A, B \subseteq \omega$, we say that $A$ is \emph{almost-contained in $B$}, written $A \subseteq^* B$, if $A \backslash B$ is finite.

\item A \emph{pseudointersection} of a family $\mathcal{F}$ of subsets of $\omega$ is a subset $P \subseteq \omega$ such that $P \subseteq^* F$ for all $F \in \mathcal{F}$.

\item A family $\mathcal{P}$ of infinite subsets of $\omega$ has the \emph{strong finite intersection property} if $\bigcap \mathcal{F}$ is infinite for all finite and non-empty $\mathcal{F} \subseteq \mathcal{P}$.
The \emph{pseudointersection number}, denoted by $\mathfrak{p}$, is the smallest cardinality of a family of subsets of $\omega$ with no infinite pseudointersection.

\item A \emph{tower} is a transfinite sequence $\transSequence{T_\alpha}{\alpha}{\lambda}$ of infinite subsets of $\omega$ such that $T_\beta \subseteq^* T_\alpha \not \subseteq^* T_\beta$ for all $\alpha < \beta < \lambda$.
The \emph{tower number}, denoted by $\mathfrak{t}$, is the smallest cardinality of a tower with no infinite pseudointersection.
\end{itemize}
\end{definition}

All three cardinals $\mathfrak{b}, \mathfrak{p}, \mathfrak{t}$ are well-defined small cardinals - see \cite{van_douwen_integers_1984}.
We will also use a `not so small' cardinal.

\begin{definition}[Novak number]
We define the \emph{Novak number}, denoted by $\mathfrak{n}$, to be the smallest cardinality of a nowhere-dense cover of $\omega^*$.
\end{definition}

We recall some facts about the Novak number from \cite{bella_sequential_2010}: $\mathfrak{t} < \mathfrak{n} \leq 2^{\mathfrak{c}}$ and it is independent of ZFC whether $2^\mathfrak{t}$ or $\mathfrak{n}$ is bounded by the other.
Moreover, every compact Hausdorff space of cardinality less than $\max(2^\mathfrak{t}, \mathfrak{n})$ is sequentially compact.

We will now use these cardinals to obtain Fr\'echet-Urysohn and sequential compactifications, provided we know that our remainder is already Fr\'echet-Urysohn and sequential respectively.
The following theorems are similar in spirit and can be summarised by the following meta-theorem: ``Any theorem that implies certain convergence properties (e.g.~sequentiality implying Fr\'echet-Urysohn, subsequentiality, sequential compactness) can be used to obtain compactification results.''

\begin{theorem} \label{Sequential compactification}
Assume $\alpha X$ is Fr\'echet-Urysohn at $\star$.
Then every compactification of $X$ with sequential remainder of cardinality strictly less than $\max(2^\mathfrak{t}, \mathfrak{n})$ is sequential.
\begin{proof}
Let $\gamma X$ be a compactification of $X$ such that $\gamma X \backslash X$ is sequential and $|\gamma X \backslash X| < 2^\mathfrak{t}$.
Note that if $A \subseteq X$ is a sequence that converges to $\star$ in $\alpha X$, then $\closure{A}^{\gamma X} = A \cup (\closure{A}^{\gamma X} \backslash X)$, so $|\closure{A}^{\gamma X}| < 2^\mathfrak{t}$ and hence is sequentially compact by \cite{bella_sequential_2010}.
Therefore by Lemma \ref{Pseudoradial compactification} $\gamma X$ is sequential.
\end{proof}
\end{theorem}

The following theorem is an adaptation of \cite[Theorem 6.2, pg.~129]{van_douwen_integers_1984}.

\begin{theorem} \label{Subsequential p}
Let $X$ be a topological space, $x \in X$ be given such that $\chi(x, X) < \mathfrak{p}$.
Then $X$ is \emph{subsequential} at $x$; that is, for all $A \in [X]^{\aleph_0}$, if $x \in \closure{A}$ then there is a sequence of $A$ that converges to $x$.
\begin{proof}
Let $A \in [X]^{\aleph_0}$ be given such that $x \in \closure{A}$.
If there exists a $y \in A \cap \nhoodCore{x}$ then $\transSequence{y}{n}{\omega} \to x$, so suppose $\nhoodCore{x} \cap A = \emptyset$.
Let $\mathcal{B}$ be a neighbourhood base for $x$ with $|\mathcal{B}| < \mathfrak{p}$.
Then $\set{A \cap B : B \in \mathcal{B}}$ has the strong finite intersection property, so there exists an infinite subset $C \subseteq A$ such that $C \subseteq^* B$ for all $B \in \mathcal{B}$.
Hence $C \to x$ and so $x$ is subsequential.
\end{proof}
\end{theorem}

\begin{theorem}
Suppose $\alpha X$ is Fr\'echet-Urysohn at $\star$ and $\chi(x, X) < \mathfrak{b}$ for all $x \in X$.
Let $\gamma X$ be a compactification of $X$ such that:
\begin{itemize}
\item $\gamma X \backslash X$ is sequential.

\item $|\gamma X \backslash X| < \max(2^\mathfrak{t}, \mathfrak{n})$.

\item $\chi(y, \gamma X) < \mathfrak{b}$ for all $y \in \gamma X \backslash X$.
\end{itemize}
Then $\gamma X$ is Fr\'echet-Urysohn.
\begin{proof}
By Theorem \ref{Sequential compactification}, $\gamma X$ is sequential, so by \cite[Proposition 3.4, pg.~534]{bella_sequential_2013} it follows that $\gamma X$ is Fr\'echet-Urysohn, since $\chi(x, X) = \chi(x, \gamma X)$ for all $x \in X$ as $X$ is locally compact and thus open in $\gamma X$.
\end{proof}
\end{theorem}

\begin{lemma} \label{Sequentially-separable compactification}
Assume $X$ is countable and let $\gamma X$ be a compactification of $X$ such that $\chi(x, \gamma X) < \mathfrak{t}$ for all $x \in \gamma X \backslash X$.
Then for all $A \subseteq X$ and $x \in \closure{A}^{\gamma X}$, there exists a sequence in $A$ that converges to $x$.
In particular, $\gamma X$ is sequentially separable\footnote{A space is \emph{sequentially separable} if it contains a countable subset $D$ such that every point is the limit of some sequence in $D$.}.
\begin{proof}
First, note that $\alpha X$ is a countable, compact Hausdorff space so is homeomorphic to an ordinal; in particular, $\alpha X$, and hence $X$, is Fr\'echet-Urysohn.
Let $A \subseteq X, x \in \closure{A}^{\gamma X}$ be given.
If $x \in X$ then $x \in \closure{A}^X$ and so there exists a sequence in $A$ that converges to $x$.
Now suppose that $x \notin X$.
By \cite{Malliaris_general_2013} and Theorem \ref{Subsequential p}, there exists a sequence in $A$ that converges to $x$.

As $X$ is dense in $\gamma X$, it follows that $\gamma X$ is sequentially separable.
\end{proof}
\end{lemma}

\begin{theorem} \label{Sequential FU compactification}
Suppose $X$ is countable and let $\gamma X$ be a compactification of $X$ such that $\gamma X \backslash X$ is Fr\'echet-Urysohn and $\chi(x, \gamma X) < \mathfrak{t}$ for all $x \in \gamma X \backslash X$.
Then $\gamma X$ is Fr\'echet-Urysohn.
\begin{proof}
Let $A \subseteq \gamma X, x \in \closure{A} \backslash A$ be given.
Notice that $\alpha X$ is a countable, compact Hausdorff space so is homeomorphic to an ordinal and hence Fr\'echet-Urysohn.
Thus $X$ is also Fr\'echet-Urysohn.

If $x \in X$ then there exists a sequence contained in $A \cap X$ that converges to $x$.
Suppose $x \notin X$.
Then $x \in \closure{A} = \closure{A \cap X} \cup \closure{A \backslash X}$.
If $x \in \closure{A \backslash X}$ then since $\gamma X \backslash X$ is Fr\'echet-Urysohn, there exists a sequence in $A$ that converges to $x$.
Otherwise, by the previous lemma there also exists a sequence in $A$ that converges to $x$.
Therefore $\gamma X$ is Fr\'echet-Urysohn.
\end{proof}
\end{theorem}

\begin{question}
Are any of the bounds in this section strict?
\end{question}

\section{Open questions}

The focus of this paper has been on building up compactifications from below, starting with the one-point compactification and extending results beyond that.
The author believes that the existence of a maximal radial / Fr\'echet-Urysohn compactication should be a fruitful line of investigation.
However, not every space has such a compactification.
To show this, we first need the following lemma:

\begin{lemma}
The infinite continuous images of $\omega + 1$ are homeomorphic to itself.
\begin{proof}
Let $f : \omega + 1 \to X$ be continuous and surjective and suppose that $X$ is infinite.
Then $X$ is homeomorphic to a infinite successor ordinal.
Without loss of generality, suppose $X = \alpha + 1$, where $\alpha$ is a countable ordinal and suppose $\alpha \geq 2 \cdot \omega$.
Then $f^{-1}[[0, \omega]]$ and $f^{-1}[[\omega + 1, 2 \cdot \omega + 1]]$ are infinite and closed, so intersect, which is a contradiction.
Thus $\alpha = \omega + n$ for some $n < \omega$ and hence $X \cong \omega + 1$.
\end{proof}
\end{lemma}

\begin{theorem}
$X \df (\omega + 1) \times \omega_1$ is a locally compact, non-compact, first-countable Hausdorff space with no maximal radial compactification, yet $\alpha X$ is radial.
\begin{proof}
First note that both $\omega + 1$ and $\omega_1$ are first-countable, so $X$ is Fr\'echet-Urysohn and hence radial.
Let $A \subseteq X$ be given such that $\star \in \closure{A}^{\alpha X}$.
Then for all $\beta < \omega_1, A \not \subseteq (\omega + 1) \times (\beta + 1)$, so there exists an $x_\beta \in A \backslash (\omega + 1) \times (\beta + 1)$.
By regularity, there exists an uncountable $B \subseteq \omega_1$ such that $\pi_{\omega + 1}(x_\beta) = \pi_{\omega + 1}(x_\gamma)$ for all $\beta, \gamma \in B$.
Now let $K \subseteq X$ be compact, so $\pi_{\omega_1}[K]$ is bounded and hence there exists a $\beta < \omega_1$ such that $K \subseteq (\omega + 1) \times (\beta + 1)$.
Then for all $\gamma \in B \backslash \beta, x_\gamma \notin K$, so $\family{x_\beta}{\beta}{B} \to \star$ and thus by Lemma \ref{Radial compactification reduction}, $\alpha X$ is radial.
Furthermore, by Corollary \ref{Finite compactification}, every finite compactification of $X$ is radial.

By \cite[Problem 3.12.20(c), pg.~237]{engelking_general_1989} the Tychonoff plank $P \df (\omega + 1) \times (\omega_1 + 1)$ is the Stone-\v{C}ech compactification of $X$ with remainder homeomorphic to $\omega + 1$, so  $X$ has no maximal finite compactifications - any compactification of $X$ is obtained by forming a closed partition of $\omega + 1$ and any finite, closed partition can be refined to a larger, finite closed partition.
Thus if $X$ has a maximal radial compactification, it must have infinite remainder.

Let $\gamma X$ be a compactification of $X$ with infinite remainder, so there exists a continuous surjection $f : P \to \gamma X$ that extends $\id_X$.
Then $f[P \backslash X] = \gamma X \backslash X$.
As $P \backslash X \cong \omega + 1$, it follows from the previous lemma that $\gamma X \backslash X \cong \omega + 1$.
Let $\delta \in \gamma X \backslash X$ be the unique non-isolated point and define $A \df \pi_{\omega + 1}[f^{-1}[\set{\delta}]], B \df X \backslash (A \times \omega_1)$.
We claim that $\delta \in \closure{B}^{\gamma X}$ but no transfinite sequence in $B$ converges to $\delta$, thus showing that there is no maximal radial compactification of $X$.

Let $U \subseteq \gamma X$ be an open neighbourhood of $\delta$.
If $f(\omega, \omega_1) \ne \delta$ then $f^{-1}[\set{f(\omega, \omega_1)}]$ is an open neighbourhood of $(\omega, \omega_1)$ and so contains $[n, \omega] \times \set{\omega_1}$ for some $n < \omega$.
But then $f[P \backslash X] = \gamma X \backslash X$ is finite, which is a contradiction.
Therefore $f(\omega, \omega_1) = \delta$ and so there exists an $n < \omega$ and an $\epsilon < \omega_1$ such that $[n, \omega] \times [\epsilon, \omega_1] \subseteq f^{-1}[U]$.
However $\gamma X \backslash X$ is infinite, so there exists an $m \in [n, \omega]$ such that $f((m, \omega_1)) \ne \delta$ and hence $(m, \epsilon) \in U$.
Thus $\delta \in \closure{B}^{\gamma X}$.

Now suppose that there exists a transfinite sequence $g$ in $B$ that converges to $\delta$.
Then $g$ converges to $\star$ in $\alpha X$, so by Lemma \ref{Ordinal remainder subsequence lemma} there exists an $m \leq \omega + 1$ and a subsequence $h$ of $g$ that converges to $(m, \omega_1)$.
However, as noted in section 2.1, no transfinite sequence in $\omega \times \omega_1$ converges to $(\omega, \omega_1)$ in $P$.
Since $f(\omega, \omega_1) = \delta$, it follows that $m < \omega$.
Moreover, $f \circ g$ converges to $\delta$, so $m \in A$.
But $\set{m} \times (\omega_1 + 1)$ is an open neighbourhood of $(m, \omega_1)$ disjoint from $B$, which is a contradiction.
Therefore no transfinite sequence in $B$ converges to $\delta$, so $\gamma X$ is not radial, concluding our proof.
\end{proof}
\end{theorem}

Regarding the structures of compactifications, the author believes the following two questions are of particular interest.

\begin{question}
When does a space have a maximal, or even greatest, radial / Fr\'echet-Urysohn compactification?
When is $\beta X$ radial / Fr\'echet-Urysohn?
\end{question}

\begin{question}
When do the radial / Fr\'echet-Urysohn compactifications form an ideal in the join-semilattice of compactifications?
\end{question}

Finally, there are still some basic questions regarding the existence of spoke systems with particular properties, most prominent is the following.

\begin{question}
When does a space have a closed, (almost-)independent\footnote{See \cite{Leek_internal_2014} for the definition of an \emph{independent} spoke system.} spoke system?
\end{question}

\bibliographystyle{alpha}
\bibliography{Bibliography}

\begin{thebibliography}{BBM13}

\bibitem[BBM13]{bella_sequential_2013}
Angelo Bella, Maddalena Bonanzinga, and Mikhail Matveev.
\newblock Sequential + separable vs sequentially separable and another
  variation on selective separability.
\newblock {\em Central European Journal of Mathematics}, 11(3):530--538, 2013.

\bibitem[BN10]{bella_sequential_2010}
Angelo Bella and Peter Nyikos.
\newblock Sequential compactness vs. countable compactness.
\newblock {\em Colloq. Math.}, 120(2):165--189, 2010.

\bibitem[Eng89]{engelking_general_1989}
Ryszard Engelking.
\newblock {\em General {T}opology}, volume~6 of {\em Sigma Series in Pure
  Mathematics}.
\newblock Heldermann Verlag, Berlin, {R}evised and completed edition, 1989.

\bibitem[Fra67]{Franklin_spaces_1967}
S.~P. Franklin.
\newblock Spaces in which sequences suffice. {II}.
\newblock {\em Fund. Math.}, 61:51--56, 1967.

\bibitem[Lee14]{Leek_internal_2014}
Robert Leek.
\newblock An internal characterisation of radiality.
\newblock {\em Topology Appl.}, 177:10--22, 2014.

\bibitem[MS13]{Malliaris_general_2013}
Maryanthe Malliaris and Saharon Shelah.
\newblock General topology meets model theory, on {$\mathfrak{p}$} and
  {$\mathfrak{t}$}.
\newblock {\em Proc. Natl. Acad. Sci. USA}, 110(33):13300--13305, 2013.

\bibitem[Sou04]{soukup_scattered_2004}
Lajos Soukup.
\newblock Scattered spaces.
\newblock In Klaas~Pieter Hart, Jun-iti Nagata, and Jerry~E. Vaughan, editors,
  {\em Encyclopedia of general topology}, pages 350--353. Elsevier Science
  Publishers, B.V., Amsterdam, 2004.

\bibitem[vD84]{van_douwen_integers_1984}
Eric~K. van Douwen.
\newblock The integers and topology.
\newblock In Kenneth Kunen and Jerry~E. Vaughan, editors, {\em Handbook of
  set-theoretic topology}, pages 111--167. North-Holland Publishing Co.,
  Amsterdam, 1984.

\end{thebibliography}

\end{document}